 \documentclass{amsart}
 \usepackage{amsmath,amsfonts,amssymb,amsthm}
 \usepackage{enumerate,color}
\usepackage[applemac]{inputenc}
 \usepackage[polutonikogreek,english]{babel}

\theoremstyle{plain}
\newtheorem{theorem}{Theorem}
\newtheorem{proposition}[theorem]{Proposition}

\newtheorem{lemma}[theorem]{Lemma}

\DeclareMathOperator{\St}{St}%
\DeclareMathOperator{\Cont}{Cont}%
\DeclareMathOperator{\diam}{diam}%

\newcommand{\tq}{\, \big| \, }%

 \title[Steiner's Invariants and Minimal Connections]{Steiner's Invariants and Minimal Connections}
 
\author[S. Ducret]{Stephen Ducret}
\address{Section de Math{\'e}matiques- IGAT,  \'Ecole Polytechnique F{\'e}derale de Lausanne, 1015 Lausanne, Switzerland}
\email{stephen.ducret@epfl.ch}

\author[M. Troyanov]{Marc Troyanov}
\address{Section de Math{\'e}matiques- IGAT,  \'Ecole Polytechnique F{\'e}derale de Lausanne, 1015 Lausanne, Switzerland}
\email{marc.troyanov@epfl.ch}

\date{June 3, 2006}

\begin{document}
\maketitle
\begin{abstract}
The aim of this note is to prove that any compact metric space can be made connected at a minimal cost, where the cost is taken to be the one-dimensional Hausdorff measure.
\end{abstract}

\section{Introduction}

Let us start to discuss the relative situation. Recall that a  \emph{continuum} is a  compact connected metric space. We denote by $\Cont(X)$ the set of all continua $C \subset X$ in an arbitrary metric space $(X,d)$

\medskip

\textbf{Definition.}  (A) Given a  metric space $(X,d)$, and a compact subset $S \subset X$, the \emph{relative Steiner invariant of $S$ in $X$} is defined as
$$
 \St(S,X)= \inf \left\{ \mathcal{H}^1(C) \tq C \mbox{ and } C \cup S  \in \Cont(X) \right\}.
$$
(B)  If $(S,d)$ is an abstract compact metric space, its \emph{absolute Steiner invariant} is defined as
\begin{multline*}
 \St(S) = \inf \left\{ \St(\iota (S),X)  \tq  X \mbox{ is an arbitrary metric space and } 
 \right.  \\  \left.\iota : S \hookrightarrow X \mbox{ is an isometric embedding} \right\}.
\end{multline*}
  
Our goal is to prove the following  
\begin{theorem}[The main theorem]
Let $S$ be a compact metric space such that $\St(S) < \infty$. Then  its Steiner invariant is realized, i.e. there exists a compact metric space $Z$ and an  isometric embedding
$S \hookrightarrow Z$ such that  $\St(\iota (S),Z)=\St(S)$. \
Furthermore, there exists $C \in \Cont(Z)$ such that $C \cup \iota (S)$ is also a continuum and $\mathcal{H}^1(C) =\St(S)$.
\end{theorem}
 
The set $C \cup \iota (S)$ described in this theorem is thus a ``minimal connection'' of $S$, i.e. it is a shortest possible set that can be added to $S$ to make it a continuum.

\section{Useful results}

It is known that the relative Steiner invariant is always realized in a proper metric space (recall that a metric space is \emph{proper} if every closed ball in it is compact):

\begin{theorem}\label{rel-th}
Let $(X,d)$ be a proper metric space, and $S \subset X$ a compact subset such that $\St(S,X) < \infty$. Then $\St(S,X)$ is realized, i.e. there exists $C \in \Cont(X)$ such that $C \cup S \in \Cont(X)$  and $\mathcal{H}^1(C) =\St(S,X)$.
\end{theorem}

\medskip

If $S$ contains only two points $\{x,y\}$, then this Theorem simply says that the two points can be joined by a shortest curve. This is the Hopf-Rinow Theorem for proper metric spaces. 

\medskip A proof of this result can be found in  \cite[Theorem 4.4.20]{ambrosio1}, it is essentially based on the Blaschke compactness Theorem for the Hausdorff distance and a semi-continuity property of the Hausdorff-measure due to  Golab. Let us recall these results.
\begin{proposition}[Blaschke]\label{blaschke}
Let $(X,d)$ be an arbitrary metric space. We denote by $\mathcal{K}(X)$ the
family of all non empty compact subsets of $X$, this is a metric space for the
Hausdorff distance $d_H$. We then have
\begin{enumerate}[a)]
  \item If $(X,d)$ is compact, then so is $(\mathcal{K}(X),d_H)$.
  \item If $(X,d)$ is proper, then so is $(\mathcal{K}(X),d_H)$.
\end{enumerate}
\end{proposition}
This theorem has been originally proved by  Blaschke in the context of convex bodies in Euclidean space. We refer to \cite[Theorem 7.3.8]{burago} or  \cite[Theorem 4.4.15]{ambrosio1} for a modern proof. 

\qed

It is not difficult to check that $\Cont(X) \subset \mathcal{K}(X)$ is a closed
subset for the topology induced by the Hausdorff distance. Furthermore:

\begin{proposition}[Golab]\label{golab}
Let $(X,d)$ be a complete metric space, and $\{C_n\} \subset \Cont(X)$ be a sequence of continua such that $C_n \rightarrow C$ for the Hausdorff distance. Then $C \in \Cont(X)$ and
\[
 \mathcal{H}^1(C) \leq \liminf_{n \to \infty}\mathcal{H}^1(C_n).
\]
\end{proposition}

See \cite[Theorem 4.4.17]{ambrosio1} for a proof. 

\qed

Our main Theorem is an extension of  Theorem \ref{rel-th}. In its proof we will need to replace the Hausdorff distance by the Gromov-Hausdorff distance and the  Blaschke Theorem will be replaced by the  Gromov compactness criterion. To recall this criterion, remember that the \emph{packing number}  of the metric space $X$ at  mesh  $\varepsilon > 0$ is the number 
\begin{multline*}
P(X,\varepsilon) = \min \left \{ n \tq \mbox{ there exists } x_1, \cdots, x_n \in X \mbox{ such that } \right. \\
\left . i \neq j \Rightarrow B(x_i,\varepsilon) \cap B(x_j,\varepsilon) = \emptyset \right \}.
\end{multline*}
Recall that metric space $X$ is totally bounded if $P(X,\varepsilon) < \infty$ for 
every $\varepsilon > 0$. 
The Gromov compactness criterion says that a family of isometry class of 
compact metric spaces is totally bounded  for the Gromov-Hausdorff distance if and only if it is \emph{uniformly totally bounded}:

\begin{theorem}[Gromov]\label{th.gromov}
Let $\mathcal{M}$ be a family of isometry classes of compact metric spaces. Then the folowing conditions are equivalent :
\begin{enumerate}[i)]
\item $\mathcal{M}$ is totally bounded for the Gromov-Hausdorff distance.
\item $\sup_{X\in \mathcal{M}} P(X,\varepsilon) < \infty$ 
for any $\varepsilon > 0$. 
\end{enumerate}
\end{theorem}

See \cite[Theorem 7.4.15]{burago}.
\qed

\bigskip

Another useful result on the Gromov-Hausdorff distance says that any sequence of compact metric spaces, which is Cauchy for the  Gromov-Hausdorff distance, contains a subsequence which can be globally realized as a subsequence of a single compact metric space:
 
\begin{proposition}\label{Embedding}
Let $\{X_n\}$ be a sequence of compact metric spaces which is a Cauchy sequence  in the Gromov-Hausdorff sense. 

Then there exists a  subsequence  $\{X_{n'}\}$, a compact metric space $Z$ and isometric embeddings $X_{n'} \hookrightarrow Y_{n'}\subset Z$ and $X \hookrightarrow Y\subset Z$ such that $Y_{n'} \rightarrow Y$ for the Hausdorff distance in $Z$.
\end{proposition}

This result is  \cite[Theorem 4.5.7]{ambrosio1}.

\qed

\section{Proof of the main Theorem}

We first need a lemma:

\begin{lemma}\label{covering_lemma}
Let $(X,d)$ be a compact metric space, and $C \in \Cont(X)$. Then for any $a \in C$ and $0 < \varepsilon < \diam(C)/2$, we have
$$
 \mathcal{H}^1(C \cap B(a,\varepsilon)) \geq  \varepsilon,
$$
in particular, we have 
$$
 P(C,\varepsilon) \leq   \frac{1}{\varepsilon} \mathcal{H}^1(C). 
$$
\end{lemma}
See  \cite[Lemma 4.4.5]{ambrosio1} for a proof of this lemma. 

\qed

We then need the following generalization of Golab's semi-continuity result: 
\begin{proposition}
Let $\{X_n\}$ be a sequence of compact metric spaces such that $X_n \rightarrow X$ in the Gromov-Hausdorff sense. Suppose that $X_n$ is connected for each $n$. Then $X$ is compact and connected, and moreover, 
$$\mathcal{H}^1(X) \leq \liminf_{n \to \infty}\mathcal{H}^1(X_n).$$
\end{proposition}

\proof  From  Proposition \ref{Embedding}, we know that, choosing  a subsequence if necessary,  there exists  a compact metric space $Z$ and isometric copies of $X_n$ and $X$ embedded in $Z$, say $Y_n, Y$, such that $Y_n\to Y$ in the Hausdorff sense. Since $X_n$ and $Y_n$ are isometric and each $Y_n$ is compact and connected, we deduce from Proposition \ref{golab} that $Y \in \Cont(Z)$ and 
$$\mathcal{H}^1(Y) \leq \liminf_{n \to \infty}\mathcal{H}^1(Y_n).$$
Now since $Y$ is isometric to $X$, we  conclude that $X$ is compact and connected as well, and that 
$$\mathcal{H}^1(X) \leq \liminf_{n \to \infty}\mathcal{H}^1(X_n),$$
because  $\mathcal{H}^1(X_n) = \mathcal{H}^1(Y_n)$ and $\mathcal{H}^1(X) = \mathcal{H}^1(Y)$.

\qed

\bigskip

\emph{Proof of Theorem 1: } 
Let $S$ be a compact metric space, and let $\{ (S_k,X_k)\}$  be a minimizing sequence for the absolute Steiner invariant of $S$, that is:
\begin{enumerate}[i.)]
\item $X_k$ is a compact metric space, and $S_k \subset X_k$ is an isometric copy of $S$;
\item $\St(S_k,X_k) \rightarrow \St(S)$.
\end{enumerate}

 \medskip

We first prove that the sequence $\left \{ X_k \right \}$ can be assumed to be uniformly totally bounded : indeed, if this where not the case, we could replace $X_k$ by $S_k \cup C_k$, where  $C_k\in \Cont(X_k)$ realizes $\St(S_k,X_k)$ (such a set exists by Theorem \ref{rel-th}).  From the compactness of $S$, we know that $P(S_k, \varepsilon) = P(S,\varepsilon) < \infty$ for all $k$ and all $\varepsilon$. Moreover, $\mathcal{H}^1(C_k) = \St(S_k,X_k)$, therefore $\sup_k \mathcal{H}^1(C_k) < \infty$ and  the family $\left \{ C_k \right \}$ is uniformly totally bounded by Lemma $\ref{covering_lemma}$. The families $\left \{ S_k \right \}$ and $\left \{ C_k \right \}$ being uniformly totally bounded, so is $\left \{ S_k \cup C_k \right \}$.

 \medskip

We henceforth assume $\left \{ X_k \right \}$   to be uniformly totally bounded.
By  the Gromov compactness criterion, Theorem \ref{th.gromov},
we know that  $\{X_k\}$ contains a subsequence which is Cauchy in the Gromov-Hausdorff distance. 
From Proposition \ref{Embedding}, we can further take a subsequence which can globally
be embedded in a compact metric space $Z$. Finally, using the Blaschke compactness Theorem, 
we can take one more subsequence, which converges for the Hausdorff distance in $Z$.
 
 \medskip
 
To sum up,  there exists a  subsequence  $\{X_{k'}\}$, a compact metric space $Z$ and isometric embeddings $\iota_{k'} : X_{k'} \hookrightarrow Y_{k'}\subset Z$ and a subset  $Y\subset Z$ such that $Y_{k'} \rightarrow Y$ for the Hausdorff distance in $Z$.

 \medskip
  
Let $T_{k'} = \iota_{k'}(S_{k'}) \subset Y_{k'}\subset Z$. By 
 Theorem \ref{rel-th}, we can find $C_{k'} \in  \Cont(Y_{k'})$ such that 
 $C_{k'}  \cup T_{k'} \in \Cont(Y_{k'} )$  and $\mathcal{H}^1(C_{k'}) =\St(T_{k'} , Y_{k'}) =
 \St(S_{k'} , X_{k'} )$.

By  Blaschke's theorem again, we may assume (taking once more a subsequence if neeed) that  $\{C_{k'}\}$ converges  to a subset $C \subset Z$ for the Hausdorff distance in $Z$.
Likewise, we may assume that  $\{T_{k'}\}$ converges  to a subset $T \subset Z$ 
(since $C_{k'}\cup T_{k'} \subset Y_{k'}$ and $Y_{k'} \to Y$, we have in fact $C\cup T \subset Y$).

Furthermore, we know from Proposition \ref{golab} that $C$ and $C\cup T$ are continua  and that
$\mathcal{H}^1(C) \leq  \liminf_{k' \to \infty}\mathcal{H}^1(C_{k'})$. But we have $\mathcal{H}^1(C_{k'}) = \St(T_{k'} , Y_{k'}) =  \St(S_{k'} , X_{k'})$, which converges to $\St(S)$, thus 
$$
\mathcal{H}^1(C) \leq  \St(S).
$$
On the other hand $C$ and $C\cup T$ are continua, hence 
$\mathcal{H}^1(C) \geq  \St(T)= \St(S)$  by definition of the Steiner invariant,
we therefore have equality.

To sum up, we have found a pair of subsets $C,T \subset Z$ such that 
$C$ and $C\cup T$ are continua, $T$ is isometric to $S$ and 
$\mathcal{H}^1(C) = \St(S)$. The proof is complete. 

\qed


\begin{thebibliography}{mm}
\bibitem{ambrosio1}  L. Ambrosio and P. Tilli \textit{Topics on Analysis in Metric Spaces}. Oxford. Lect. Series. Math. 25.  (2004)

\bibitem{papadopoulos} A. Papadopoulos,
 \textit{Metric Spaces, Convexity and Positive Curvature}. 
IRMA Lectures in Mathematics and Theoretical Physics.  (2005)
 

\bibitem{burago} D. Burago, Y. Burago, S. Ivanov
\textit{A Course in Metric Geometry.}
Graduate studies in Mathematics. American Mathematical Society (2001)

 \end{thebibliography}
\end{document}